\documentstyle[12pt, amstex]{amsart}

\renewcommand{\pf}{\noindent{\em Proof: }}
\def\s#1{\cal{#1}}

      \newtheorem{Thm}{Theorem}[section]
\newtheorem{Def}[Thm]{Definition} \newtheorem{Rem}[Thm]{Remark}
\newtheorem{Lem}[Thm]{Lemma} \newtheorem{Cor}[Thm]{Corollary}
\newtheorem{Prop}[Thm]{Proposition}

\newcommand{\ra}{\rangle}
\newcommand{\la}{\langle}
\newcommand{\Span}{\operatorname{Span}}
\newcommand{\g}[1]{\mathfrak{#1}}

\begin{document}
\author{Hans Plesner Jakobsen}
\title{Q-differential operators}
\date{\today}
\maketitle

\begin{abstract}
We set up a framework for discussing ``$q$-analogues'' of the usual covariant
differential operators for hermitian symmetric spaces. This turns out
to be directly related to the deformation quantization associated to  quadratic
algebras satisfying certain conditions introduced by Procesi and De Concini.
\end{abstract}

\section{Introduction}

The investigation, of which we are here reporting some results, began
with the question about what should be ``quantized wave operators'' in
the context of (quantized) hermitian symmetric spaces. Immediately,
there is a very simple thing one can do, namely one can pass to the
quantized enveloping algebra. Here, there are unitarizable highest
weight modules and the most singular of these have kernels which, in
analogy with the case $q=1$ can be said to be ``quantized wave
operators'' (\cite{MR97c:81074}, \cite{jakobsen;unitarity-quantum}).
However, when $q$ is generic, there is no immediate space of functions
on which these differential operators act.

The first objects we have come across in our attempt to repair on this
are (families of) quadratic algebras that seem to replace the
hermitian symmetric spaces. See \cite{MR1667369},
\cite{jakobsen;quantized-hermitian-symmetric},  and below. Secondly, a natural
setting for differential operators (in an algebraic approach) could be
duality. Combining these two one comes across the following:

\medskip

Let ${\s P}$ be a projection (not necessarily self adjoint)
in the tensor algebra $T(V)$ over some (finite-dimensional) vector
space. Suppose that ${\s P}$ maps $T^r(V)$ to $T^r(V)$ for each $r$.
Solutions ${\mathcal P}$ of ($\star$) or ($\star\star$) to the following
equations, reminiscent of the Yang-Baxter equations, turn out to have a
fundamental importance.
\[
\boxed{
\begin{array}{rcl}(\star)\quad\forall r,s: (I_r\otimes{\s  P}\otimes I_s){\s P}&=&{\s P}\\
 (\star\star)\quad\forall r,s:{\s P}(I_r\otimes{\s P}\otimes
I_s)&=&{\s P}\end{array} }\] Indeed, such a partial solution can be
used to define an associative algebra of polynomial functions on
either $V$ (case of ($\star$)) or $V^*$ (case of ($\star\star$)). And,
once this has been established, one may introduce, by duality,
quantized differential operators.

\medskip

We will construct below, for a quadratic algebra that satisfies a
certain technical condition, a projection $\s P$ (a quantized
symmetrization map) which solves both equations at the same time. The
condition is related to the condition in ``the Diamond Lemma'' by
Bergman (\cite{bergman;diamond-lemma})  -- a major influence for us in
relation to this part. The condition is satisfied by the quadratic
algebras connected with hermitian symmetric spaces.

One aspect of some of the quadratic algebras that fulfill the
condition (including those from hermitian symmetric spaces) is that
they give rise to Poisson structures. The deformed products obtained
from $\s P$ is directly related to this in the usual way.   

\medskip

Our way of quantizing holomorphic functions may be extended to all
functions by quantizing anti-holomorphic functions independently and
then, based on considerations involving e.g. reproducing kernels,
letting holomorphic and anti-holomorphic variables commute. We mention 
that other possibilities have been extensively studied by, in
particular, D. Shklyarov, S. Sinel'shchikov, and L. Vaksman. See e.g. \cite{shklyarov-vaksman}
and references cited therein, or math.QA. at 
http://xxx.lanl.gov/.
\medskip

The material is organized as follows: In Section~\ref{prev} we give a
short description of the way covariant differential operators arise in
the classical case, c.f.
\cite{harris-jakobsen;covariant-differential-operators},
\cite{jakobsen;basic-covariant-differential-operators}. In
Section~\ref{q-alg} and Section~\ref{q-alg-t},  quadratic algebras are 
introduced, examples are given, and some technical assumptions are
discussed. Then, in Section~\ref{q-sym} the operators $\mathcal P$ are 
finally constructed and basic properties are given. The associative
(non-commutative) polynomial algebras are introduced via duality in
Section~\ref{dual}, and it is briefly discussed how different choices
of bases may give different presentations of the same algebra. The
quantized differential operators are then  introduced by means of
the duality. In Section~\ref{mq2} the situation is analyzed in detail
for $M_q(2)$. The ``$q$-differential operators'' are seen to consist
of some rather agreeable components together with possibly a more
complicated term which points towards covariant differentiation in
infinite dimensional spaces. Further aspects of this will be presented 
in forthcoming papers. Finally, some computations of the differential 
operators for $M_q(n)$ are appended.

\medskip

\section{The classical situation or how to get the wave operator, the
  Dirac operator, Maxwell's equations etc. (in the mass 0 case/absence
  of sources case) without physics.}

\label{prev}
Let $\s B$ be an irreducible hermitian symmetric space of the
noncompact type. Then (c.f. Helgason \cite[Chapter VIII]{MR26:2986} ${\s B}$ is
diffeomorphic to $G/K$ where $G$ is a connected noncompact simple Lie
group with trivial center and $K$ is a maximal compact subgroup with
non-discrete center. If $\g g,\g k$ denote the complexified Lie
algebras of $G,K$, respectively, then there are complex subalgebras
$\g p^\pm$ such that
\begin{eqnarray}
\g g &=& \g p^-\oplus \g k\oplus \g p^+,
\\\nonumber
[\g p^\pm,\g p^\pm]&=&0,\\\nonumber
[\g p^+,\g p^-]\subseteq \g k,&\textrm{ and }& [\g k, p^\pm]\subseteq \g p^\pm.
\end{eqnarray}
Moreover, we choose a subalgebra ${\g h}$ which is both a Cartan
subalgebra for $\g g$ and $\g k$.
Observe that we have
\begin{equation}
\boxed{\s U(\g g)= \s U(\g p^-)\cdot \s U(\g k)\cdot\s U(\g p^+).}
\end{equation}

Operationally, it is here more adequate to use the equivalent
description where $\s B$ is a bounded symmetric domain in ${\mathbb
  C}^N$, $G$ is the connected component of the group of biholomorphic
bijections of $\s B$ onto itself, and $K$ is the isotropy group of a
point. Indeed, we may, and shall, take $\s B$ to be an open bounded
subset of $\g p^-$ such that $0\in \s B$ and such that $K$ acts
linearly.

\medskip

Let $\tau$ be a unitary representation of $K$ in a finite dimensional
vector space $V_\tau$. Then $G\times_K V_\tau$ is a vector bundle over
$\s B$ and $G$ acts naturally on the space
$\Gamma_{\textrm{h}}(G\times_K V_\tau)$ of holomorphic sections of
$G\times_K V_\tau$. The bundle is equivalent to a trivial bundle ${\s
B}\times V_\tau$ and as a result one obtains a representation $U_\tau$
of $G$ in the space ${\s H}({\g p}^-)\otimes V_\tau$ of $V_\tau$
valued holomorphic functions on $\s B$. The algebraic span of the $K$
types is exactly the space ${\s P}({\g p}^-)\otimes V_\tau$ of
$V_\tau$ valued polynomials on $\g p^-$. 

\medskip

We say that a differential operator \begin{equation*} {\s D}: {\s
    H}({\g p}^-)\otimes V_{\tau_1}\mapsto{\s H}({\g p}^-)\otimes
  V_{\tau_2}\end{equation*}
\noindent is covariant provided\begin{equation}
  \forall g\in G:U_{\tau_2}(g){\s D}={\s D}U_{\tau_1}(g).
\end{equation}
It follows from the assumptions that $\s D$ is a holomorphic constant
coefficient $\hom(V_{\tau_1},V_{\tau_2})$ valued differential
operator.

It turns out that such operators indeed do exist, even under
additional unitarity assumptions, but to get a better understanding of
where they come from, we turn to another construction:

\begin{Def}For $V_\tau$ as before,
\[M(V_\tau)={\s U}({\g g})\otimes_{{\s U}(\g k +\g p^+)}V_\tau\]
is called a generalized Verma module. It is a highest weight module
generated by a non-zero highest weight vector $v_{\tau}$.
Specifically, ${\g p}^+v_\tau={\g k}^+v_\tau=0$ and $\forall h\in{\g
  h}:h\cdot v_\tau=\Lambda_\tau(h)\cdot v_\tau$ for some (highest
weight) $\Lambda_\tau\in{\g h}^*$.
\end{Def}

The analogue of a covariant operator at this level is a ${\s U}(\g g)$
homomorphism $\phi:M(V_{\tau_2})\mapsto M(V_{\tau_1})$. A homomorphism
$\phi$ is completely determined by $\hat v_{\tau_2}=\phi(1)$ - a
vector in $M(V_{\tau_1})$ which has the same weight as $v_{\tau_2}$
and which is annihilated by ${\g p}^+$ and ${\g k}^+$. Conversely any
such so called primitive vector (for physicists: a secondary vacuum)
determines a homomorphism.

\medskip

The key fact now is the following

\begin{Prop}
There is a natural non-degenerate pairing between ${\s P}\otimes
V_\tau$ and $M(V_{\tau^\prime})=M(V_{\tau}^\prime)$ under which the
spaces as ${\s U}(\g g)$ modules are the dual to each other. Under
this duality, a homomorphism between generalized Verma modules
correspond to a covariant differential operator in the dual picture -
and conversely, a covariant differential operator determines in the
dual picture a homomorphism.
\end{Prop}

\medskip

Another key fact is that there occur naturally some homomorphisms
between generalized Verma modules at singular unitary holomorphic
representations. Indeed, the homomorphism is defined in terms of the
lowest ``missing $\g k$ type.

\medskip

\medskip

Consider the symmetric algebras 
\begin{equation}
\boxed{S({\g p}^\pm)=T({\g p}^\pm)/I_\pm(XY-YX)}
\end{equation}
where $I_\pm(XY-YX)$ denotes the ideal in $T({\g p}^\pm)$ generated by
all elements of the form $X\otimes Y-Y\otimes X$ with $X,Y\in {\g
p}^\pm$. This is clearly a quadratic algebra. Let ${\s P}_0^\pm$
denote the projections of $T({\g p}^\pm)$ onto $S({\g p}^\pm)$. These
are well known maps:
\begin{equation}
\boxed{{\s P}_0^\pm\textrm{ are symmetrization maps.}}
\end{equation}

The Killing form $B$ on ${\g g}$ gives a non-degenerate pairing between
${\g p}^+$ and ${\g p}^-$. For $w^+\in {\g p}^+$ and $z^-\in {\g p}^-$
we write $\la w^+,z^-\ra=B(w^+,z^-)$. This extends to a pairing
between $S({\g p}^+)$ and $S({\g p}^-)$ by 
\begin{eqnarray}
\la [w^+_1\otimes\cdots\otimes w_r^+], [z^-_1\otimes\cdots\otimes
z_s^-]\ra=\delta_{r,s}\sum_{\sigma\in S_r}\prod_{i=1}^r\la
w_i^+,z^-_{\sigma(i)}\ra\\\nonumber =\delta_{r,s}r!\la
w^+_1\otimes\cdots\otimes w_r^+, {\s P}_0^-(z^-_1\otimes\cdots\otimes
z_s^-)\ra.
\end{eqnarray}
Through this pairing, any element $[w]=[w^+_1\otimes\cdots\otimes
w_r^+]\in S({\g p}^+)$ defines a polynomial ${\s F}^0_{[w]}\in {\s
  P}({\g p}^-)$ by
\begin{equation}
\boxed{{\s F}^0_{[w]}(z^-)=\la {\s P}^+_0 (w^+_1\otimes\cdots\otimes w_r^+),
z^-\otimes\cdots z^-\otimes\cdots\ra.}
\end{equation}
In this way we get an identification of vector spaces (indeed, of ${\g k}$ modules)
\begin{equation}
{\s P}({\g p}^-)\otimes V_\tau=S({\g p}^+)\otimes V_\tau.
\end{equation}
Similarly, 
\begin{equation}
M(V_\tau^\prime)=S({\g p}^-)\otimes V_\tau^\prime,
\end{equation}
and the pairing between the two modules is just the introduced pairing
between $S({\g p}^+)$ and $S({\g p}^-)$ augmented with the pairing
between the module $V_\tau$ and its dual module $V_\tau^\prime$.

\medskip
 
Now observe that ${\g p}^-$ acts on ${\s P}({\g p}^-)$ by contraction,
\begin{equation}
\boxed{\begin{array}{ccl}(z_0^-{\s F}^0_{[w]})((z^-))&=&\la (w^+_1\otimes\cdots\otimes
w_r^+), {\s P}^-_0(z_0^-\otimes z^-\otimes\cdots
z^-\otimes\cdots\ra \\ &=& (\frac{\partial}{\partial z_0^-}{\s
F}^0_{[w]})((z^-)).\end{array}}
\end{equation} 

But this is just a differentiation, and in this way, $S({\g p}^-)$ can
be viewed as either a space of polynomials on ${\g p}^+$ or as a space
of constant coefficient differential operators on ${\g
p}^-$. Extending the above to the case of generalized Verma modules,
$M(V_\tau^\prime)$ is the space of $V_\tau^\prime$ valued constant
coefficient differential operators on ${\g p}^-$ and the pairing above
can be formulated as follows: If $z^-\mapsto p^-(z^-)\otimes v\in {\s
P}({\g p}^-)\otimes V_\tau$ and if $p^+(\frac{\partial}{\partial
z^-})\otimes v^\prime\in M(V_\tau^\prime)$ then
\begin{equation}
\la p^-\otimes v,p^+\otimes
v^\prime\ra=\left(p^+\left(\frac{\partial}{\partial
z^-}\right)p^-\right)(0)\cdot\la v,v^\prime\ra.
\end{equation}

\medskip

Finally, observe that the product in e.g. ${\s P}({\g p}^-)$ is given by
\begin{equation}
\boxed{{\s F}^0_{[w^a]}\star{\s F}^0_{[w^b]}={\s F}^0_{[w^a\otimes w^b]}.}
\end{equation}
The well-definedness of this follows from
\begin{equation}\label{0cru1}
\boxed{(I_r\otimes{\s P}_0^+\otimes I_s){\s P}_0^+={\s P}_0^+(I_r\otimes{\s
    P}_0^+\otimes I_s)={\s P}_0^+.}
\end{equation}

This star is commutative simply because we work with real symmetrization.

\medskip

\section{Quadratic algebras}
\label{q-alg}
Our construction below, though inspired by hermitian symmetric spaces,
works for a more general class of algebras, namely quadratic algebras
(subject to some technical assumptions to be stated later). We first
give some examples and then later the precise definitions.

\subsection{Examples of quadratic algebras}The simplest quadratic
algebras are the commutative ones 
\[X_iX_j-X_jX_i=0.\] We see that the symmetric algebras of the
previous section fall in this category. Other examples are ``quantized 
objects'' e.g. 

\[\begin{array}{rcll}
  AB&=&qBA\qquad &\textrm{ (quantum plane)}\\
AB-q^2BA&=&1\qquad &\textrm{ (quantized Weyl)},
\end{array}\]
where $q\in{\mathbb C}^*$ is the quantum parameter.

One of the most studied ones is the quantized function algebra of
$n\times n$ matrices, $M_q(n)$, defined by the relations
\begin{eqnarray}\label{aiii}{\textbf A}{\mathbf III:} \qquad Z_{i,j}Z_{i,k}&=&qZ_{i,k}Z_{i,j} \text{ if } j<k,\\\nonumber
Z_{i,j}Z_{k,j}&=&qZ_{k,j}Z_{i,j} \textrm{ if }i<k,\\\nonumber
Z_{i,j}Z_{s,t}&=&Z_{s,t}Z_{i,j} \textrm{ if }i<s\textrm{ and }t<j,\\\nonumber
Z_{i,j}Z_{s,t}&=&Z_{s,t}Z_{i,j}+(q-q^{-1})Z_{i,t}Z_{s,j} \text{ if }
i<s\textrm{ and }j<t.
\end{eqnarray}

This algebra is in the class of quadratic algebras connected with
quantized hermitian symmetric spaces and for this reason we sometimes refer to
it as  {\bf AIII}. We mention two more from the class, namely  {\bf CI} (but
{\bf DIII} is also covered by this) and {\bf BDI}($q=2$). Observe
that a misprint in the relations for   {\bf CI} has been corrected and two missing
relations have been added compared to (\cite{jakobsen;quantized-hermitian-symmetric})
\begin{eqnarray*}{\textbf C}{\mathbf  I:}\qquad\quad\qquad\quad\\
W_{i,i}W_{j,j}-W_{j,j}W_{i,i} &=&\frac{1-q^2}{q+q^{-1}} W_{i,j}^2\text{ (}i<j\text{)},
\\
W_{i,i}W_{j,k}-W_{j,k}W_{i,i} &=&(1-q^2)W_{i,j}W_{i,k}\text{ (}i<j\text{ and 
}j<k\text{)}, \\
qW_{i,j}W_{j,k}-W_{j,k}W_{i,j} &=&(q^{-2}-q^2)W_{i,k}W_{j,j}\text{ (}i<j<k
\text{)}, \\
W_{i,j}W_{k,l}-W_{k,l}W_{i,j} &=&q^{-1}W_{i,k}W_{j,l}-qW_{i,k}W_{j,l}\text{ (
}i<j<k<l\text{)}, \\
W_{i,j}W_{k,k}-W_{k,k}W_{i,j}&=&(1-q^2)W_{i,k}W_{j,k}, \\
W_{i,j}W_{k,l}-W_{k,l}W_{i,j} &=&(q^{-1}-q)W_{i,l}W_{k,j}\text{ (}i,k<j<l
\text{)}, \\
W_{i,i}W_{i,j} &=&q^{-2}W_{i,j}W_{i,i}\text{ (}i<j\text{)}, \\
W_{i,,j}W_{j,j} &=&q^{-2}W_{j,j}W_{i,j}\text{ (}i<j\text{)}, \\
W_{i,j}W_{i,k} &=&q^{-1}W_{i,k}W_{i,j}\text{ (}i<j<k\text{)}, \\
W_{i,k}W_{j,k} &=&q^{-1}W_{j,k}W_{i,k}\text{ (}i<j<k\text{)}, \\
W_{i,j}W_{k,l} &=&W_{k,l}W_{i,j}\text{ (}i\leq j,k<i,\text{ and }j<l\text{)}.
\end{eqnarray*}

The relations for type {\bf BDI} are: 
\begin{eqnarray}
{\textbf B}{\textbf D} {\mathbf I:} \qquad W_iW_{i+r+1}&=&q^{-1}W_{i+r+1}W_i\text{ if }r\geq 0\text{ and }r\neq 2(n-i),
\quad\quad\\
W_iW_{2n+1-i}-W_{2n+1-i}W_i&=&-qW_{i+1}W_{2n-i}+q^{-1}W_{2n-i}W_{i+1}\\&& \qquad\quad \qquad\quad\qquad \text{
for }1=1,\ldots ,n-1.
\end{eqnarray}

If we let $Z_i-1=(-q)^{-i}W_i$ for $i=1,\dots,n$ and $Z_i^*=W_{2n-i}$
for $i=0,\dots,n-1$ the last relations are seen (replacing $n$ by $N$) to
be those of the quantized Heisenberg space (by some called the quantum
symplectic space - a name which according to our classes is somewhat
confusing), $F_q(N)$ of the quantum space ${\mathbb C}^N$, i.e. the
associative algebra generated by $z_0, z_1,\cdots,z_{N-1},z_0^*
z_1^*,\cdots,z_{N-1}^*$ subject to the following relations:
\begin{eqnarray}\label{jzrem} z_iz_j&=&q^{-1}z_jz_i\text{ for }i<j,\\\nonumber
      z_i^*z_j^*&=&qz_j^*z_i^*\text{ for }i<j,\\\nonumber
 z_iz_j^*&=&q^{-1}z_j^*z_i\text{ for }i\ne j,\textrm{ and }\\\nonumber
 z_iz_i^*-z_i^*z_i&=&(q^2-1)\sum_{k>i}z_kz_k^*.\end{eqnarray}

\medskip

\medskip

\subsection{General definition}

The general definition of a quadratic algebra is as follows:

Let $V$ denote an $N$-dimensional complex vector space, let $T=T(V)$
denote the tensor algebra over $V$, let $R$ be a subspace of $V\otimes
V$, and let $I_R$ denote the 2-sided ideal in $T$ generated by the
$R$. Then
\begin{Def}
\[{\mathcal A}= T/I_R.\]
\end{Def}
\noindent We say that $I_R$ is the space generated by the relations.

\medskip

The starting point of our present investigation is the following fact: 

Let $\g g$, $\g k$ be as in the Section~\ref{prev}. Let ${\s U}_q(\g g)$ and ${\s U}_q(\g k)$
be the
quantized enveloping algebras of $\g g$ and  $\g k$,
respectively. Then there are quadratic algebras ${\s A}^\pm$ which 
furthermore are ${\s U}_q(\g k)$ modules such that

\begin{equation}
\boxed{{\s U}_q(\g g)= {\s A}^-\cdot {\s U}_q(\g k)\cdot{\s A}^+.}
\end{equation}

The quadratic algebras satisfy the additional assumptions below.

\medskip

\section{Technical discussion}
\label{q-alg-t}
We consider a quadratic algebra ${\mathcal A}$ generated by (linearly 
independent) elements
$X_1,\dots,X_N$.  For each $i=1,\dots,N$ let ${\mathcal A}_i$ denote
the algebra generated by $X_1,\dots,X_i$. We assume that
the defining relations are of the form:
\begin{equation}\nonumber\label{straight}
\text{(Rel)\qquad \quad If }i>j\text{ then }X_{i}X_{j}=b_{ij}X_{j}X_{i}+p_{ij}, \text{ with
  }p_{ij}\in {\mathcal A}_{i-1}.\end{equation}

Let $V$ denote the $N$-dimensional complex vector space
spanned by the elements $X_{1}, \ldots,X_{N}$, let $T=T(V)$ denote the tensor algebra
over $V$, and let $I_R$ denote the ideal in $T$ generated by elements $
X_{i}X_{j}-(b_{ij}X_{j}X_{i}+p_{ij})$. Then 
\begin{equation}
{\mathcal A}:=T/I_R.
\end{equation}

\medskip

For $r\in \Bbb{N}$  we let $T^{r}=\underbrace{V\otimes \cdots \otimes
  }_{r}V$. To an element $X=X_{i_1}\otimes X_{i_2}\otimes\cdots\otimes
X_{i_r}\in T^r$ we associate the element $\ell(X)=(n_1,n_2,\dots,n_N)\in \{0,1,\dots,r\}^N$ where 
\begin{equation}\forall i=1,\dots,N: n_i=\#\{s\mid i_s=i\}.
\end{equation}

\medskip

We shall from now on drop the $\otimes$ whenever this can be done
without placing the presentation in jeopardy.

\medskip

We now introduce lexicographic ordering $\leq_l$ on
$\{0,1,\dots,r\}^N$ (according to which $(0,\dots,0,r)$ is the biggest
element) to introduce a partial ordering, also denoted $\leq_l$, on
the set of monomials in $T^r$ simply by declaring $u^\prime\leq_l
u\Leftrightarrow \ell(u^\prime)\leq_l\ell(u)$.

\medskip

The essential assumption (EA) which we now make is introduced to avoid
situations where, due to some special cancelations, a sum of elements
in $I_R$ might add up to an element which strictly precedes all the
summands in the order. Specifically, we assume for any element $u\in
I_R$
\begin{equation}\nonumber\begin{array}{l}
\text{(EA): }\\\label{ea}
u\in\Span\{a\cdot(X_iX_j-b_{ij}X_jX_i-p_{ij})\cdot b\mid a,b \in
T\text{ and }a\cdot(X_iX_j)\cdot b\leq_l u\}.
\end{array}
\end{equation}

\medskip

In the following we shall introduce certain operations which are
related to thinking of (Rel) as a reduction system. The reductions are
then of the form
\begin{equation}
 (X_iX_j,b_{ij}X_jX_i +p_{ij})\quad(\textrm{ for all }i>j).
\end{equation}
Indeed, we can, analogously to \cite[Section
3]{bergman;diamond-lemma}, introduce a {\em misordering index} $i(Z)$
of an element $Z=X_{i_1}\dots X_{i_r}$ as the number of pairs of
indices $(i_a,i_b)$ in $Z$ for which $i_a>i_b$. This we can combine
with the ordering $\leq_l$ to give a new partial ordering, $\leq$ on
monomials as follows:
\begin{equation}
u_1<
u_2\stackrel{\textrm{Def.}}{\Leftrightarrow}\left\{\begin{array}{l}u_1<_lu_2
\textrm{ or}\\\ell(u_1)=\ell(u_2)\textrm{ and }i(u_1)<i(u_2).\end{array}\right.
\end{equation}

It is clear that if $A,B$ are monomials in $T$, then
$u^\prime<u\Rightarrow A\cdot u\cdot B<A\cdot u^\prime\cdot B$. Thus,
our partial ordering is a {\em semigroup partial ordering}. Moreover,
for all $r,s$ with $s>r$ we have that $b_{sr}X_rX_s $ is of strictly
smaller misordering index and $p_{sr}$ is of strictly smaller lexicographic order
than $X_sX_r$. Thus $b_{sr}X_rX_s+ p_{sr}$ is of strictly less order
(w.r.t $<$) than $X_sX_r$ and hence,  the reduction system is {\em
  compatible with the reduction system}.

The two mentioned properties are parts of the requirements for the
Diamond Lemma \cite[Theorem 1.2]{bergman;diamond-lemma} to be applicable to our
situation.

\begin{Prop}
All elements of $T$ are reduction unique.
\end{Prop} 

\pf By observing that all reductions decrease the order it follows
that the system satisfies the descending chain condition.  It remains,
according to \cite[p. 181]{bergman;diamond-lemma}, to prove that all
ambiguities of the reduction system are resolvable. The only place where we can get
ambiguities are on terms $X_iX_jX_k$  with $i>j>k$. Here we must
prove (still following \cite[p. 181]{bergman;diamond-lemma}) 
\begin{equation}
Y=(b_{ij}X_jX_i+p_{ij})X_k-X_i(b_{jk}X_kX_j+p_{jk})\in
I_{i,j,k},
\end{equation}
where $I_{i,j,k}$ denotes the subspace spanned by all elements
$A((X_sX_r-b_{sr}X_rX_s -p_{sr})B$ with $s>r$ and
$A(X_sX_r)B<X_iX_jX_k$. But clearly,
$(b_{ij}X_jX_i+p_{ij})X_k-X_i(b_{jk}X_kX_j+p_{jk})\in I_R$ (it is the
reduction of $X_iX_jX_k- X_iX_jX_k$). Secondly, the only monomials in
$Y$ that map to $\ell(X_iX_jX_k)$ under the map $\ell$ are
$b_{ij}X_jX_iX_k$ and $b_{jk}X_iX_kX_j$. But after two more
reductions, they both become $b_{ij}b_{jk}b_{ik}X_kX_jX_i$ plus
something of lower lexicographic order. Since the two original terms
have opposite signs, the highest order terms cancel. The claim then
follows from (EA). According to the Diamond Lemma we are done. \qed

\medskip

We immediately get

\begin{Cor}
The set $\{X_1^{i_1}\cdots X_N^{i_N}\mid 1_1,\dots, i_n\in {\mathbb N}_0\}$ is
a basis for ${\mathcal A}$.
\end{Cor}

\begin{Cor}
${\mathcal A}$ is a domain and is in fact an iterated twisted
polynomial algebra. In
particular, the assumptions of Procesi and De Concini
(\cite{MR95j:17012})  are satisfied.
\end{Cor}

\medskip

Conversely we have the following result which implies that the
algebras {\bf AIII, BDI}, and  {\bf CI} above fit into the framework:

\begin{Prop}
Given a quadratic algebra ${\mathcal A}$ as above, satisfying (Rel), and  furthermore satisfying

\smallskip

\textrm{$($DCP$)$}\quad In all cases where $j<i$ set
$\sigma_i(X_j)=b_{ij}X_j$. Then for each $i$, $\sigma_i$ defines an
automorphism of ${\mathcal A}_{i-1}$.

\smallskip

Then it  satisfies \textrm{$($EA$)$}.
\end{Prop}

\pf As in (\cite{MR95j:17012}) it follows that the algebra is an iterated twisted
polynomial algebra. Suppose that (EA) is not satisfied. Let $u\in I_R$
be the smallest element which does not satisfy (EA). Then up to this
order, the algebra behaves exactly as an iterated twisted polynomial
algebra. But the advent of $u$ then implies that there is at least one
extra relation at this level. But this contradicts the fact that the
algebra has the same Hilbert series as its associated quasipolynomial
algebra (the algebra where the relations are $X_iX_j=b_{ij}X_jX_i$).
\qed

\medskip

\begin{Rem}
  It would be interesting to classify all quadratic algebras that
  satisfy this reduction assumption (EA) or, equivalently, (DCP). It is clearly a quite strong assumption, on the
  order of complication of e.g. the Jacobi Identity in the enveloping
  algebra. 
\end{Rem}

\medskip

In \cite[Theorem 1.2]{bergman;diamond-lemma}, Bergman goes on to
define a product and projection etc. but we are after something else -
though also a projection.
\label{tech}
\medskip

\section{The construction}
\label{q-sym}

\medskip

Maintain the notation of Section~\ref{tech}.

\begin{Def}\label{sdef} We define a linear map $S:V\otimes V\longrightarrow V\otimes V$ by 
\begin{eqnarray}
S(X_{i}\otimes X_{j}) &=&b_{ij}X_{j}X_{i}+p_{ij}\text{ if }i>j,
\\\label{quasi}
S(X_{j}\otimes X_{i}) &=&(b_{ij})^{-1}(X_{i}X_{j}-p_{ij})\text{ if
  }i>j, \textrm{ and}  \\
S(X_{i}\otimes X_{i}) &=&X_{i}\otimes X_{i}\text{ for all }i=1,\ldots ,N.
\end{eqnarray}
Furthermore, we define $\overline{S}:V\otimes V\longrightarrow V\otimes V$
by 
\begin{eqnarray}
\overline{S}(X_{i}\otimes X_{j}) &=&b_{ij}X_{j}X\text{ if }i>j, \\
\overline{S}(X_{j}\otimes X_{i}) &=&(b_{ij})^{-1}(X_{i}X_{j})\text{ if
  }i>j, \textrm{ and} 
\\
\overline{S}(X_{i}\otimes X_{i}) &=&X_{i}\otimes X_{i}\text{ for all }
i=1,\ldots ,N.
\end{eqnarray}
\end{Def}

From now on, we assume that 
\begin{equation}
\forall i,j:b_{ij}=q^{\alpha _{ij}},
\end{equation}
where $q$ until further notice is a non-zero complex number. Recall that the
associated quasi-polynomial algebra is the quadratic algebra $\overline{
{\mathcal A}}$, generated (for clarity) by elements $x_{1,}\ldots ,x_{N\text{ 
}}$ with relations $x_{i}x_{j}=q^{\alpha _{ij}}x_{j}x_{i}$.

\begin{Def}
For $i\in {\mathbb N}, \sigma_{i}$ denotes the linear map $T\longrightarrow T$ given by 
\begin{eqnarray}
&&\sigma_{i}(v_{1}\otimes \ldots \otimes v_{i-1}\otimes v_{i}\otimes
v_{i+1}\otimes \cdots \otimes v_{n} \\
&=&v_{1}\otimes \ldots \otimes v_{i-1}\otimes S(v_{i}\otimes v_{i+1})\otimes
\cdots \otimes v_{n}
\end{eqnarray}
and $\overline{\sigma}_{i}$ denotes the linear map $T\longrightarrow T$
given by 
\begin{eqnarray}
&&\overline{\sigma }_{i}(v_{1}\otimes \ldots \otimes v_{i-1}\otimes
v_{i}\otimes v_{i+1}\otimes \cdots \otimes v_{n} \\
&=&v_{1}\otimes \ldots \otimes v_{i-1}\otimes \overline{S}(v_{i}\otimes
v_{i+1})\otimes \cdots \otimes v_{n}\text{.}
\end{eqnarray}
\end{Def}

\smallskip

We now state and prove a series of lemmas about these maps.  

\begin{Lem}
For each $i$,
\begin{equation}
u_1\leq_l u_2\Leftrightarrow \overline{\sigma
  }_{i}(u_1)\leq_l \overline{\sigma }_{i}(u_2)\Leftrightarrow {\sigma
  }_{i}(u_1)\leq_l {\sigma }_{i}(u_2).  
\end{equation}
\end{Lem}

\pf Clear from the definitions. \qed

\medskip

\begin{Lem}
\label{lem1}For each $i\in \Bbb{N}$, $\sigma _{i}$ is equal to the identity
modulo $I_R$, i.e. for each $u\in T$ there exists an $r\in I_R$ such that 
\begin{equation}
\sigma _{i}(u)=u+r\text{.}
\end{equation}
\end{Lem}

\noindent{\em Proof:} This is obvious from the definitions. \qed

\smallskip

\begin{Lem}
\label{lem2}For each $i\in \Bbb{N}$, $\overline{\sigma }_{i}$ $\overline{
\sigma }_{i+1}$ $\overline{\sigma }_{i}=$ $\overline{\sigma }_{i+1}$ $
\overline{\sigma }_{i}$ $\overline{\sigma }_{i+1}$, and hence $\overline{
\sigma }_{1},\ldots ,\overline{\sigma }_{n-1}$ define a representation,
called {\bf quasi-permutation}, of the symmetric group S$_{n}$ on $T^{n}$.
\end{Lem}

\pf
By choosing $a_{ij}\,$appropriately, we may write $\overline{S}(X_{i}\otimes
X_{j})=q^{a_{ij}}X_{j}\otimes X_{i}$ for all $i,j=1,\ldots ,N$. The claim
follows easily from this by an elementary computation. \qed

\medskip
From now, in all statements involving {\em order}, we mean the lexicographical order $\leq_l$.
\medskip

\begin{Lem}
$\label{lem3}$For each $i\in \Bbb{N}$, $\sigma _{i}=\overline{\sigma }_{i}$
modulo lower order.
\end{Lem}

\pf
Obvious from the definitions. \qed

\begin{Lem}
\label{lem4}The following hold

\begin{enumerate}
\item  For each $i\in \Bbb{N}$, $\sigma _{i}$ preserves $I_R$.

\item  For each $i\in \Bbb{N}$, if for $u\in T:\overline{\sigma }_{i}(u)=u$,
then $\sigma _{i}(u)=u$.

\item  $\sigma _{i}\sigma _{i+1}\sigma _{i}=\sigma _{i+1}\sigma _{i}\sigma
_{i+1}$ modulo $I_R\,$or modulo lower order terms.
\end{enumerate}
\end{Lem}

\smallskip 

\pf
The first claim follows from Lemma~\ref{lem1}. To prove the
second claim it is clearly enough to prove that if for $u\in V\otimes V\,$, $
\overline{S}(u)=u$, then $S(u)=u$, and for this, we may assume that 
\begin{equation}
u=X_{i}\otimes X_{j}+q^{a_{ij}}X_{j}\otimes X_{i}\text{.}
\end{equation}
The assertion then follows by an easy computation. The validity of the
part of the last statement that involves $I_R$ follows from
Lemma~\ref{lem1} combined with the first item of this lemma. The
validity of the other part follows from Lemma~\ref{lem3} combined with
Lemma~\ref{lem2}. \qed

\begin{Lem}
\label{lem5}Let $u=a\cdot(X_{i}X_{k}-b_{ik}X_{k}X_{i}-p_{ik})\cdot b
\in I_R$, where $a,b\in{\mathcal A}$ and $a$ is a homogeneous polynomial of degree
$j-1$. Then there exists a positive integer $p$ such that 
$(1+\sigma_j )^{p}u=0$. 
\end{Lem}

\smallskip 

\pf
We have that 
\begin{equation}
(1+\sigma_j)\left(a\cdot(X_{i}X_{k}-b_{ik}X_{k}X_i-p_{ik})\cdot
b\right)=a\cdot(1-\sigma_j)p_{ik}\cdot b\text{.}
\end{equation}
Since clearly, by construction and by Lemma~\ref{lem4},
$(1-\sigma_j)p_{ik}$ is of lower order, is in $I_R$, and is of the right
form ((EA) is not needed here) one may repeat the procedure with $u$
replaced by $u^\prime=a\cdot(1-\sigma_j)p_{ik}\cdot b$. After a finite
number of steps one will reach $0$. \qed

\medskip

We now wish to introduce an analogue of the usual symmetrization map on
$T$.  Let us first consider the representation of $S_{n}\,$described
in Lemma~\ref{lem2}. For any $\sigma \in S_{n}\,$we denote the
resulting operator on $T^{n}$ as $\overline{\sigma }$ and we set
\begin{equation}
P_{\text{quasi-sym}}=\frac{1}{n!}\sum_{\sigma \in S_{n}}\overline{\sigma }\text{,
}
\end{equation}
and call this operator {\bf quasi-symmetrization}. It is clear that this
operator is the projection onto the subspace of tensors in $T^{n}\,$that are
invariant under each $\overline{\sigma }_{i},i=1,\ldots ,n-1$. More
precisely, the following identities of course hold just as for ordinary
symmetrization:

\begin{Lem} \label{5.9}
\[\forall i: \qquad\overline{\sigma }_i\cdot P_{\text{quasi-sym}}=P_{\text{quasi-sym}}\cdot
\overline{\sigma }_i=P_{\text{quasi-sym}} . 
\]
\end{Lem}

\medskip

We next want to define a similar operator on $T^{n}\,$with respect to the $
\sigma _{i}$'s. The problem is, of course, that we do not have a bona fide
representation. In spite of this we proceed by defining for each $\sigma \in
S_{n}$ an operator $\widehat{\sigma }=\sigma _{i_{1}}\sigma _{i_{2}}\cdots
\sigma _{i_{r}}$ if $\sigma =s_{i_{1}}s_{i_{2}}\cdots s_{i_{r}}$, where $
s_{j},j=1,\ldots ,n-1$, denotes the elementary transpositions in $S_{n}$ and
we set, for each such set of decompositions of elements, 
\begin{equation}
P=\frac{1}{n!}\sum_{\sigma \in S_{n}}\widehat{\sigma }.
\end{equation}

Notice that for each $i=1,\ldots ,n-1$ we have a left coset decomposition
of $S_{n}$ with respect to the subgroup $\left\{ 1,s_{i}\right\}
$; $S_{n}=C_{i}\times \left\{ 1,s_{i}\right\} $ for some suitable
subset $C_i$ of $S_{n}$. Hence we have, among
the operators $P$, some of the form (all denoted $P_{i}$)
\begin{equation}
P_{i}=\widetilde{P}\cdot (1+\sigma _{i}).
\end{equation}
More generally we can introduce 
\begin{equation}
P_{i,r}=\widetilde{P}\cdot \left(\frac{1+\sigma _{i}}2\right)^{r}
\end{equation}
where $r$ later will be taken to be a sufficiently big power.

\smallskip

\begin{Cor}
\label{cor1}Each $P$ leaves $I_R$ invariant and $P=P_{\text{quasi-sym}}$ modulo
lower order.
\end{Cor}

\smallskip

\pf
This follows directly from Lemma~\ref{lem3} and Lemma~\ref{lem4}. \qed

\medskip

\begin{Lem}
\label{lem6}Let $u\in I_R$. Then there exists an $N\in \Bbb{N}$ such that 
\begin{equation}
P^{N}(u)=0.
\end{equation}
\end{Lem}

\pf By linearity and by Lemma~\ref{lem5}, we may assume that
$P_{i,r}(u)=0$ for some $i,r$. But then, since $P$ and $P_{i}$ agree
modulo lower order, $P(u)$ is of lower order that $u$. And by
Corollary~\ref{cor1} $P(u)\in I_R$. Now invoke (EA) to yield that we
after finitely many steps  get $ P^{N}(u)=0$. \qed

\medskip

\begin{Cor}\label{cor2.2}
  If $u\in I_R$ satisfies  $P_{\text{quasi-sym}}(u)=u$ then $u=0$.
\end{Cor}

\pf Combine Lemma~\ref{lem4} (2.) with Lemma~\ref{lem6}. \qed

\medskip

\begin{Lem}
\label{lem7}Let $u\in T$. Then there exists an $N_{0}\in \Bbb{N}$ such that 
\begin{equation}
P^{N}(u)=\widetilde{P}^{N}(u)
\end{equation}
for all $N\geq N_{0}$.
\end{Lem}

\smallskip

\pf We have that $P(u)=P_{\text{quasi-sym}}(u)+u_{1}$ where $u_{1}$ is of
lower order.  By Lemma~\ref{lem3} and Lemma~\ref{lem4} it follows that
$P^{2}(u)=P_{q\text{-sym }}(u)+P_{\text{quasi-sym}}(u_{1})+u_{2}.$ Thus,
there exists a $\widehat{u}$ such that
$P^{N}(u)=P_{\text{quasi-sym}}({\widehat{u}})$. Likewise, there exists a $
\widehat{u}$ such that
$\widetilde{P}^{N}(u)=P_{\text{quasi-sym}}(\widetilde{u})$.  Moreover,
clearly
\begin{equation}
P_{\text{quasi-sym}}(\widetilde{u}) =P_{\text{quasi-sym}}(\widehat{u})\mod I_R,
\end{equation}
and hence, by Corollary~\ref{cor2.2} the claim follows. \qed

\medskip

\begin{Def}
  Set
\begin{equation}
{\mathcal P}_{q\text{-sym}}=\lim_{N\longrightarrow \infty }P^{N}.
\end{equation}
\end{Def}

The following is immediate

\smallskip

\begin{Prop}
\label{prop1}${\mathcal P}_{q\text{-sym}}$ is a well-defined projection satisfying
\begin{equation}
{\mathcal P}_{q\text{-sym}}(I_R)=0.
\end{equation}
\end{Prop}

\smallskip

\begin{Lem}
\label{lem8}If $P(u)=0$, then $u\in I_R$. If ${\mathcal P}_{q\text{-sym}}(u)=0$ then $u\in I_R$.
\end{Lem}

\smallskip

\pf This follows directly from Lemma~\ref{lem1}. \qed

\begin{Lem}
\label{cor2}If ${\mathcal P}_{q\text{-sym}}(u_{1})=0$ then ${\mathcal P}_{q\text{-sym}}(u_{1}\otimes u)=0\,$
for all $u\in T$.
\end{Lem}

\smallskip

\pf It follows by Lemma~\ref{lem8} that $u_{1}\in I_R$. Hence
$u_{1}\otimes u\in I_R$ . The claim then follows from
Proposition~\ref{prop1}. \qed

\medskip

We shall occasionally denote the restriction of ${\mathcal P}_{q\text{-sym}}$ to $T^k$ by
${\mathcal P}_{q\text{-sym}}^k$, but most of the times we drop the subscript. For $r,s,k\in
{\mathbb N}$ define the linear operator $I_r\otimes {\mathcal P}_{q\text{-sym}}^k\otimes
I_s$ from $T$ into $T$ by

\begin{eqnarray}I_r\otimes {\mathcal P}_{q\text{-sym}}^k\otimes I_s(v_1\otimes\cdots\otimes v_r\otimes
  v_{r+1}\otimes\dots v_{r+k}\otimes v_{r+k+1}\cdots\otimes v_{r+k+s})\\\nonumber
v_1\otimes\cdots\otimes v_r\otimes{\mathcal P}_{q\text{-sym}}^k(
  v_{r+1}\otimes\dots v_{r+k})\otimes v_{r+k+1}\cdots\otimes v_{r+k+s}
\end{eqnarray}

\medskip

The crucial property of ${\mathcal P}_{q\text{-sym}}$ then is

\begin{Prop}
\begin{equation}\label{cru1}
\boxed{(\star)\quad \forall r,s: (I_r\otimes {\mathcal P}_{q\text{-sym}}^k\otimes I_s){\mathcal P}_{q\text{-sym}}^{r+k-s}={\mathcal P}_{q\text{-sym}}^{r+k-s}.}
\end{equation}
\end{Prop}

\pf As in the proof of Lemma~\ref{lem7} observe that for any $u\in T$,
${\mathcal P}_{q\text{-sym}}(u)$ is quasi-symmetric. Hence the claim follows directly from
Lemma~\ref{lem4} and Lemma~\ref{5.9}. \qed

\smallskip

We also have

\begin{Prop}
\begin{equation}\label{cru2}
\boxed{(\star\star)\quad \forall r,s:{\mathcal P}_{q\text{-sym}}^{r+k-s}(I_r\otimes {\mathcal P}_{q\text{-sym}}^k\otimes I_s)={\mathcal P}_{q\text{-sym}}^{r+k-s}.}
\end{equation}
\end{Prop}

\pf By Proposition~\ref{prop1}, it suffices to prove
that for any $u\in T$, $(I_r\otimes{\mathcal P}_{q\text{-sym}}^k\otimes I_s)(u)-u\in {I_R}$. Here, it
suffices to consider a $u$ of the form $u_1\otimes \dots\otimes
u_r\otimes v\otimes v_1\otimes \dots\otimes v_{s}$ with $v\in T$. Then
\begin{equation}
(I_r\otimes{\mathcal P}_{q\text{-sym}}^k\otimes I_s)(u)-u=u_1\otimes\dots\otimes u_r\otimes
({\mathcal P}_{q\text{-sym}}^k(v)-v)\otimes v_1\otimes \dots\otimes v_{s}
\end{equation} 
and the claim follows from Lemma~\ref{lem8} since by construction,
${I_R}$ is an ideal in $T$.  \qed

\medskip

\medskip

\begin{Rem}\label{sa-proj}
It is of course possible to introduce an inner product in $T(V)$ in
which the projection ${\mathcal P}_{q\text{-sym}}$ is self-adjoint. Indeed, there is an
infinite family of possible choices. It remains to be decided, if there is  a {\em natural} candidate.
\end{Rem}

\medskip

\section{Duality}
\label{dual}

\subsection{New observations}

We maintain the assumptions on $\s A$. Let $V^*$ denote the linear dual to $V$ and denote the pairing by
\begin{equation}
V^*\times V\ni v^*,v\mapsto\la v^*,v\ra.
\end{equation}
We extend this pairing to a pairing between $T^*=T(V^*)$ and $T$ in
the usual tensor product fashion.

Clearly, the introduced structure can be transported to $T^*$ by this
duality. On the level of the pairing between $V\otimes V$ and
$V^*\otimes V^*$, we can consider the transposed of the $S$ and
$\overline S$ of Definition~\ref{sdef}. More generally, we can
consider the projection $({\mathcal P}_{q\text{-sym}})^t$ on $T^*$. Let $I_R^t$ denote the
kernel of the restriction of $({\mathcal P}_{q\text{-sym}})^t$ to $V^*\otimes V^*$ and use
$I_R^t$ to define a quadratic algebra ${\s A}^t$. 

\begin{Prop}\label{quasiprop} ${\s A}^t$ is a quasipolynomial algebra.  
\end{Prop}

\pf This follows from condition (\ref{quasi}) which implies that the
columns in the matrix of ${\mathcal P}_{q\text{-sym}}$ corresponding to $X_i\otimes
X_j$ and $b_{ij}\cdot X_j\otimes X_i$ have simple sums and
differences. The transposed then have the same property for rows and
this immediately gives that any pair $X_i^*,X_j^*$ satisfies a
quasipolynomial identity. Of course, there might a priori be more
relations than that, but this is ruled out by dimension considerations
in the dual algebra. \qed

\begin{Rem}Proposition~\ref{quasiprop} is perhaps surprising to the
  point of being disappointing. Notice however that the result relies
  on the chosen duality between $T(V)$ and $T(V^*)$. Other choices,
  e.g. based on inner products  as in Remark~\ref{sa-proj} combined
  with a conjugation,  may perhaps lead
  to other algebras, but  we shall not pursue this point here. 
\end{Rem}

\medskip

For $w\in {T^n}^*$ and $z\in T^n$ we
define the   ${\mathcal P}_{q\text{-sym}}$-symmetrized pairing $\la\la\cdot,\cdot\ra\ra$ by
\begin{equation}
\la\la w,z\ra\ra=n!\la w,{\mathcal P}_{q\text{-sym}}(z)\ra .
\end{equation}
This is the pairing that generalizes the pairing
$(q,p)=(q(\frac{\partial}{\partial z}),(p(\cdot))(0)$ between
polynomials and differential operators.

\begin{Def}\label{abs-def}For $w\in {T^*}$ and $z\in T$:
\begin{equation}
\boxed{{\mathcal F}_w(z)=\la w, {\mathcal P}_{q\text{-sym}}(z)\ra.}
\end{equation}
\end{Def}
It is clear that
\begin{equation}
{\mathcal F}_w(z)={\mathcal F}_{[w]}([z])
\end{equation}
where $[z]$ and $[w]$ denote the equivalence classes in ${\mathcal A}$  and
${\mathcal A}^t$, respectively, corresponding to $z$ and $w$.

\medskip

\begin{Def}\label{def}
\begin{equation}
\boxed{({\mathcal F}_{[w_1]}\star {\mathcal F}_{[w_2]}):={\mathcal
  F}_{[w_1\otimes w_2]}.}
\end{equation}
\end{Def}

\medskip

\begin{Prop}
  The product in Definition~\ref{def} is a well-defined associative
  product.
\end{Prop}

\pf The associativity is clear as soon as it is well-defined. This it
is by (\ref{cru1}). \qed

\medskip

\begin{Rem}
Of course, there is the expected direct relation between the Poisson structure
defined by the above non-commutative product,

$$\lim_{q\rightarrow 1}\frac1{q-1} \left({\mathcal F}_{[w_1]}\star
  {\mathcal F}_{[w_2]}- {\mathcal F}_{[w_1]}\star {\mathcal
    F}_{[w_2]}\right),$$ and the usual Poisson structure for certain
quadratic algebras as defined by Procesi and De Concini
(\cite[p. 84-85]{MR95j:17012}).
\end{Rem}

\medskip

We consider ways of representing  the functions ${\mathcal F}_{[w]}$
as functions on $V$. 

Let $X_1,\dots,X_N$ be a basis of $V$ as in (Rel) in Section~\ref{tech}
and let 
\begin{equation}
\forall \alpha=(\alpha_1,\dots,\alpha_N)\in({\mathbb
  N}_0)^N: X^\alpha=X_1^{\alpha_1}\otimes\dots\otimes X_N^{\alpha_N}. \end{equation} 
Furthermore, choose for each $\alpha\in({\mathbb
  N}_0)^N$ a homogeneous polynomial 
\[p^{\s C}_\alpha=c_\alpha z_1^{\alpha_1}\dots
z_N^{\alpha_N}+\sum_{\beta<\alpha}d_{\alpha,\beta} z^\beta\] where each
$d_{\alpha,\beta}$ is a complex number, where each $c_\alpha$ is a
non-zero constant, and where the symbol $\s C$ (e.g. a lower triangular $\infty\times
\infty$ matrix) represent these choices. The ordering  $\beta<\alpha$ is lexicographic.

\begin{Def}\label{c-def}
\[\boxed{{\s F}^{\s C}_{[w]}(z_1,\dots z_N):=\la w,{\mathcal P}_{q\text{-sym}}(\sum_\alpha 
p^{\s C}_\alpha X^\alpha)\ra.}\]
\end{Def}

The following is immediate

\begin{Prop}
For each ${\s C}$ we get a faithful representation of the algebra ${\mathcal
  A}^t$ in an associative algebra of polynomial functions on ${\mathbb C}^N$.
\end{Prop}

\smallskip

The family of algebras we have defined by means of $\s C$ includes
algebras defined by other (PBW-like) bases since a change of basis
will simply be equivalent to a change of $\s C$. For some specific choice of
$\s C$'s, a given element may give rise to a differential operator of an 
especially simple form, c.f. Section~\ref{mq2} below.

\medskip

\begin{Rem}We shall only pursue certain specific versions of
  Definition~\ref{c-def} below, but we wish to mention here that one
  may in fact go even further and represent the abstract functions of
  Definition~\ref{abs-def} as non-commutative functions with values in
  certain algebras. In doing this, the construction is related to some
  algebras occurring when $q$ is an $m$th root of unity. Specifically, for
$M(n,{\mathbb C})$,   observe that
\begin{equation}\label{poly-alg}
X_{1,1}^{a_{1,1}}\cdots X_{n,n}^{a_{n,n}}X_{1,1}^{m\cdot
  b_{1,1}}\cdots X_{n,n}^{m\cdot b_{n,n}} 
\end{equation}
with $0\leq a_{i,j}\leq m-1$ for all $1\leq i,j\leq n$ form a
basis of $\s A$ for each $m\in{\mathbb N}$. Suppose namely that we could write 0
as a non-trivial linear combination of these. The coefficient of the
highest order term is then by definition non-zero. However, we can
rewrite the basis elements with respect to the standard basis. Doing
this, the highest order term remains unchanged. But then the
coefficient must be zero since the other basis is indeed a basis. Thus
the elements are linearly independent, and by considering degrees,
they must be a spanning set. 

We can then interpret a specific  element  $X_{1,1}^{c_{1,1}}\cdots X_{n,n}^{c_{n,n}}X_{1,1}^{m\cdot
  d_{1,1}}\cdots X_{n,n}^{m\cdot d_{n,n}} $ in (\ref{poly-alg}) as corresponding 
to the polynomial $z_{1,1}^{d_{1,1}}\cdots z_{n,n}^{d_{n,n}}\otimes
(X_{1,1}^{c_{1,1}}\cdots X_{n,n}^{c_{n,n}})$ with values in the space spanned by the elements
$X_{1,1}^{a_{1,1}}\cdots X_{n,n}^{a_{n,n}}$ with $0\leq a_{i,j}\leq m-1$ for all $1\leq i,j\leq n$.
\end{Rem}

\medskip

We now consider,  for
$M(n,{\mathbb C})$, some specific instances of Definition~\ref{c-def}:
\begin{Def}If $\{z_{i,j}\}_{i,j=1}^n\in M(n,{\mathbb C})$ set
  $z=\sum_{i,j=1}^n z_{i,j}X_{i,j}$ and
\begin{eqnarray}{\mathcal F}_{[w]}^{(1)}(z_{11},\dots,z_{n,n})&=&\la w,{\mathcal
    P}_{q\text{-sym}}(Z\otimes \cdots \otimes Z)\ra\\\nonumber
{\mathcal F}_{[w]}^{(2)}(z_{11},\dots,z_{n,n})&=&\la w,{\mathcal
    P}_{q\text{-sym}}(\sum_\alpha c_\alpha z^\alpha X^\alpha)\ra, 
\end{eqnarray} where 
\begin{equation}
c_\alpha=\frac{(\vert\alpha\vert)!}{(\alpha_{1,1})!\cdot\cdots\cdot(\alpha_{n,n})!}.
\end{equation}
\end{Def}

\medskip

Now, let $[w_\beta]$ be determined by 
\begin{equation}{\mathcal F}_{[w_\beta]}^{(2)}(z_{1,1},\dots,z_{n,n})=z^\beta,
\end{equation}
i.e.
\begin{equation}\la ({\mathcal P}_{q\text{-sym}})^t(w_\beta), X^\alpha\ra =(c_\beta)^{-1}\delta_{\alpha,\beta}.
\end{equation}
By duality we have
\begin{equation}
\boxed{\frac{\partial}{\partial X_0}{\mathcal
  F}^{(i)}_{[w_\beta]}(\cdot)=\frac1{(\vert\beta\vert -1)!}\la\la w_\beta,X_0(\cdot)\ra\ra.}  
\end{equation}
Thus, 
       
\begin{equation}\left(\frac{\partial}{\partial X_0}{\mathcal
      F}^{(i)}_{[w_\beta]}\right)(Z)=
    \left \{\begin{array}{l}\vert\beta\vert\cdot\la w_\beta,{\mathcal P}_{q\text{-sym}}(X_0\otimes Z\otimes\cdots\otimes 
Z)\ra\text{ for } i=1 \\\vert\beta\vert\cdot\la
      w_\beta,{\mathcal P}_{q\text{-sym}}(X_0\otimes(\sum_\alpha c_\alpha z^\alpha X^\alpha))\ra \text{ for }i=2
    \end{array} \right. .
\end{equation}

If our ordering of $X$ is $X_{1,1},X_{1,2},\dots,X_{n,n}$ then we get in
particular that 
\begin{equation}
\frac{\partial}{\partial X_{1,1}}{\mathcal
  F}^{(2)}_{[w_\beta]}(Z)=\vert\beta\vert\frac{c_\alpha z^\alpha}{c_\beta}\delta_{\alpha-1,\beta}=\beta_{1,1}z_{1,1}^{\beta_{1,1}-1}z_{1,2}^{\beta_{1,2}}\cdot z_{n,n}^{\beta_{n,n}}. 
\end{equation}

Likewise, 

\begin{eqnarray}
\frac{\partial}{\partial X_{1,2}}{\mathcal
  F}^{(2)}_{[w_\beta]}(Z)&=&\vert\beta\vert\frac{c_\alpha z^\alpha}{c_\beta}\delta_{\alpha-1,\beta}=q^{-\beta_{1,1}}\beta_{1,2}z_{1,1}^{\beta_{1,1}}z_{1,2}^{\beta_{1,2}-1}\cdot z_{n,n}^{\beta_{n,n}},\\\nonumber \frac{\partial}{\partial X_{2,1}}{\mathcal
  F}^{(2)}_{[w_\beta]}(Z)&=&\vert\beta\vert\frac{c_\alpha
  z^\alpha}{c_\beta}\delta_{\alpha-1,\beta}=q^{-\beta_{1,1}}\beta_{2,1}z_{1,1}^{\beta_{1,1}}z_{1,2}^{\beta_{1,2}}z_{2,1}^{\beta_{2,1}-1}\cdot z_{n,n}^{\beta_{n,n}},\\\nonumber\text{and (for $2\times2$) case}\\\nonumber\frac{\partial}{\partial X_{2,2}}{\mathcal  F}^{(2)}_{[w_\beta]}(Z)&=&\beta_{2,2}q^{(-\beta_2-\beta_3)}z_{1,1}^{\beta_{1,1}}z_{1,2}^{\beta_{1,2}}z_{2,1}^{\beta_{2,1}}z_{2,2}^{\beta_{2,2}-1}\\\nonumber &-q& (1-q^{-2 \beta_{1,1}-2 })\frac{\beta_{1,2}\beta_{2,1}}{\beta_{1,1}+1} z_{1,1}^{\beta_{1,1}+1}z_{1,2}^{\beta_{1,2}-1}z_{2,1}^{\beta_{2,1}-1}z_{2,2}^{\beta_{2,2}}.
\end{eqnarray}

\section{$M_q(2)$}\label{mq2}

We continue with the functions ${\mathcal F}_{[w_\beta]}^{(2)}(z_{1,1},\dots,z_{n,n})$ from the previous section but specialize 
further to the quantized function algebra of $2\times2$ matrices.

\medskip

Let $z_1=z_{1,1}, z_2=z_{1,2}, z_3=z_{2,1}$, and $z_4=z_{2,2}$. Then
\begin{eqnarray*}
\left(\frac{\partial}{\partial z_1}\right)_q
(z_1^{\alpha_1}z_2^{\alpha_2}z_3^{\alpha_3}z_4^{\alpha_4})&=&\alpha_1z_1^{\alpha_1-1}z_2^{\alpha_2}z_3^{\alpha_3}z_4^{\alpha_4}\\\left(\frac{\partial}{\partial
    z_2}\right)_q
(z_1^{\alpha_1}z_2^{\alpha_2}z_3^{\alpha_3}z_4^{\alpha_4})&=&q^{-\alpha_1}\alpha_2z_1^{\alpha_1}z_2^{\alpha_2-1}z_3^{\alpha_3}z_4^{\alpha_4}\\
\left(\frac{\partial}{\partial z_3}\right)_q
(z_1^{\alpha_1}z_2^{\alpha_2}z_3^{\alpha_3}z_4^{\alpha_4})&=&q^{-\alpha_1}\alpha_3z_1^{\alpha_1}z_2^{\alpha_2}z_3^{\alpha_3-1}z_4^{\alpha_4}\\
\left(\frac{\partial}{\partial z_4}\right)_q
(z_1^{\alpha_1}z_2^{\alpha_2}z_3^{\alpha_3}z_4^{\alpha_4})&=&q^{-\alpha_2-\alpha_3}\alpha_4z_1^{\alpha_1}z_2^{\alpha_2}z_3^{\alpha_3}z_4^{\alpha_4-1}\\
&+&\alpha_2\alpha_3{\s
  K}_{\alpha_1}q^{-2\alpha_1+2}z_1^{\alpha_1}z_2^{\alpha_2-1}z_3^{\alpha_3-1}z_4^{\alpha_4}
\end{eqnarray*}
where, with $q=e^\hbar$, \begin{eqnarray}{\s
    K}_{\alpha_1}&=&-q^{2\alpha_1-1}(1-q^{-2\alpha_1-2})\frac{z_1}{\alpha_1+1}\\\nonumber &=&-e^{-3\hbar}(2\hbar+\dots+\frac{(2\hbar)^n}{n!}(\alpha+1)^{n-1}+\cdots)\cdot
  z_1.\end{eqnarray} If we let $S_1= z_1\frac{\partial}{\partial z_1}$ then we
see that ${\s  K}_{\alpha_1}={\s K}_1$ independently of $\alpha_1$ where the operator 
\begin{equation} \label{43}{\s
  K}_1\equiv -e^{-3\hbar}(2\hbar+\dots+\frac{(2\hbar)^n}{n!}(S_1)^{n-1}+\cdots)\cdot
z_1\end{equation}only involves the variable $z_1$. 
The factors $q^{-\alpha_1}$ and $q^{-\alpha_2-\alpha_3}$ can of course also be
dealt with analogously. However, if we define
\begin{equation}\label{44}K_i(z_1^{\alpha_1}z_2^{\alpha_2}z_3^{\alpha_3}z_4^{\alpha_4})=q^{-\alpha_i}z_1^{\alpha_1}z_2^{\alpha_2}z_3^{\alpha_3}z_4^{\alpha_4}\quad
  i=1,2,3,4\end{equation}
then these are just like the usual $K$ operators and we may then also
write (if there is no subscript on a differential operator it means that it is a
classical differential operator)
\begin{eqnarray}\label{4-dif}
\left(\frac{\partial}{\partial z_1}\right)_q&=&\left(\frac{\partial}{\partial
    z_1}\right) \\\nonumber
\left(\frac{\partial}{\partial z_2}\right)_q&=&K_1\left(\frac{\partial}{\partial
    z_2}\right)\\\nonumber 
\left(\frac{\partial}{\partial
    z_3}\right)_q&=&K_1\left(\frac{\partial}{\partial
    z_3}\right) \\\nonumber
\left(\frac{\partial}{\partial
    z_4}\right)_q&=&K_2K_3\left(\frac{\partial}{\partial z_4}\right)+{\s K}_1
 \left (\frac{\partial}{\partial z_2}\right)_q\left(\frac{\partial}{\partial z_3}\right)_q\\\nonumber
\left(\frac{\partial}{\partial
    z_4}\right)_q&=&K_2K_3\left(\frac{\partial}{\partial z_4}\right)+{\s O}_1
 \left (\frac{\partial}{\partial z_2}\right)\left(\frac{\partial}{\partial z_3}\right),
\end{eqnarray}
where ${\s O}_1={\s K}_1K_1^2$. 

Notice that $e^{-\hbar S}=K_1$, $\frac{\partial}{\partial z_1}\cdot {\s
K}_1=-e^{-3\hbar}\left(e^{\hbar(2S+2)}-1\right)=(q^{-3}-q^{-1}K_1^{-2})$, 
and   $\frac{\partial}{\partial z_1}\cdot{\s
O}_1=(q^{-3}K_1^{2}-q^{-1})$.

The operators $\left(\frac{\partial}{\partial z_i}\right)_q,
i=1,2,3,4$, satisfy similar relations as  (\ref{aiii}) for $X_{1,1},
X_{2,1},X_{1,2,},X_{2,2}$ except that $q\rightarrow q^{-1}$. In
particular, what corresponds to the wave operator $\square_q$ is the central
element 

\begin{equation}
\square_q=\left(\frac{\partial}{\partial
    z_1}\right)_q\left(\frac{\partial}{\partial z_4}\right)_q-q^{-1}\left(\frac{\partial}{\partial z_2}\right)_q\left(\frac{\partial}{\partial z_3}\right)_q.
\end{equation}

It is perhaps somewhat surprising that in this case the mixed degrees 
disappear again and
\begin{equation}
\square_q=K_2K_3\frac{\partial}{\partial z_1}\frac{\partial}{\partial z_4}-q\frac{\partial}{\partial z_2}\frac{\partial}{\partial z_3}.
\end{equation}

However, in the case of e.g. the Dirac operator, which basically will
be a $2\times2$ matrix with entries (up to constant multiples)
$\left(\frac{\partial}{\partial z_1}\right)_q,\dots,
\left(\frac{\partial}{\partial z_4}\right)_q$, there is no
cancelation of the second order term arising from $\left(\frac{\partial}{\partial z_4}\right)_q$.

\medskip

We now discuss further the first order differential operators of
(\ref{4-dif}). First of all we remark that the simple appearance of
$\left(\frac{\partial}{\partial z_1}\right)_q$ is a result of the
given choice of ordering. Other choices of orderings
(or, equivalently, of constants $c_\alpha$) can make the other
variables have a simple appearance - at the expense of that of $z_1$. 

Secondly introduce the coordinate functions 
\begin{equation}{\mathcal
    F}_{i}^{(2)}(z_{1},\dots,z_{4})=z_i\textrm{ for }i=1,2,3,4.
\end{equation}

We can then introduce the left derivatives $\delta^L_i$ for
$i=1,2,3,4$:

\begin{eqnarray}
&\left(\delta^L_i {\mathcal
    F}_{[w]}^{(2)}\right)(z_1,\dots,z_4)=&\\\nonumber&\lim_{u\rightarrow 0}\left(\left({\mathcal    F}_{i}^{(2)}(ue_i)\right)^{-1}\star\left({\mathcal
    F}_{[w]}^{(2)}((z_1,\dots,z_4)+ue_i)-{\mathcal
    F}_{[w]}^{(2)}(z_1,\dots,z_4)\right)\right)&
\end{eqnarray}

It follows easily that $\left(\frac{\partial}{\partial z_i}\right)_q=\delta^L_i$ for
$i=1,2,3$ and $\left(\frac{\partial}{\partial
    z_4}\right)_q=\delta^L_4+{\mathcal K}_1 \delta^L_2\delta^L_3$.

\medskip

We finish this section with a study of how in particular  $\left(\frac{\partial}{\partial z_4}\right)_q$
may be viewed as a covariant derivative. Let ${\hat{\s O}}={\s O}_1\left
  (\frac{\partial}{\partial z_2}\right)\left(\frac{\partial}{\partial
    z_3}\right)$. Set
\begin{equation}F(f)=\left(\begin{array}{c}f\\{\hat{\s O}}f\\({\hat{\s
          O}}^2+ [{\hat{\s O}},K_2K_3\left(\frac{\partial}{\partial
    z_4}\right)])f\\({\hat{\s
          O}}^3+[{\hat{\s
          O}}^2,K_2K_3\left(\frac{\partial}{\partial
    z_4}\right)]+K_2K_3\left(\frac{\partial}{\partial
    z_4}\right)[{\hat{\s O}},K_2K_3\left(\frac{\partial}{\partial
    z_4}\right)])f\\\vdots
\end{array}\right)
\end{equation}
Let \begin{equation}A=\left(\begin{array}{ccccc}0&1&0&0&\dots\\0&0&1&0&\dots\\\vdots&\vdots&\vdots&\vdots&\vdots\end{array}\right)\end{equation}

Then \begin{equation}(K_2K_3\left(\frac{\partial}{\partial
    z_4}\right)+A)F(f)=F(\left(\frac{\partial}{\partial
    z_4}\right)_qf).\end{equation}

\medskip

Another possibility is to let 
\begin{equation}G(f)=\left(\begin{array}{c}f\\K_4^2\left
  (\frac{\partial}{\partial z_2}\frac{\partial}{\partial
    z_3}\right)f\\K_4^4\left
  (\frac{\partial}{\partial z_2}\frac{\partial}{\partial
    z_3}\right)^2f\\K_4^6\left
  (\frac{\partial}{\partial z_2}\frac{\partial}{\partial
    z_3}\right)^3f\\\vdots
\end{array}\right)
\end{equation}
and 
\begin{equation}B=\left(\begin{array}{cccccc}0&{\s O}_1K_4^{-2}&0&0&0&\dots\\0&0&{\s O}_1K_4^{-2}&0&0&\dots\\0&0&0&{\s O}_1K_4^{-2}&0&\dots\\\vdots&\vdots&\vdots&\vdots&\vdots\end{array}\right).\end{equation}

\medskip

Then we also have

\begin{equation}(K_2K_3\left(\frac{\partial}{\partial
    z_4}\right)+B)G(f)=G(\left(\frac{\partial}{\partial
    z_4}\right)_qf).\end{equation}

This last version is also well behaved with respect to the other generators.

\medskip

\section*{Appendix}

Here we compute the first order differential operators for  $M_q(n)$. Observe that we have

\begin{eqnarray}\label{difform}
&&z_{n,i}(z_{a,1}^{\alpha_{a,1}}z_{a,2}^{\alpha_{a,2}}\cdots
z_{a,n}^{\alpha_{a,n}})=\\\nonumber
&&\sum_{x=1}^{i-1}c_{a,x}q^{(\alpha_{a,x+1}+\cdots +\alpha_{a,i-1})}z_{a,1}^{\alpha_{a,1}}\cdots z_{a,x}^{\alpha_{a,x}-1}\cdots
z_{a,i}^{\alpha_{a,i}+1}\cdots z_{a,n}^{\alpha_{a,n}}z_{n,x}\\\nonumber&&+q^{\alpha_{a,i}}(z_{a,1}^{\alpha_{a,1}}z_{a,2}^{\alpha_{a,2}}\cdots
z_{a,n}^{\alpha_{a,n}})z_{n,i},
\end{eqnarray}
where $c_{a,x}=q(q^{-2\alpha_{a,x}}-1)$ and where the exponent to $q$
should be interpreted as 0 for $x=i-1$.

With (\ref{difform}) to our disposal we can now give the general form
of $\frac{\partial}{\partial X_{i,j}}$. Let
$\Gamma_{i,j}=\Gamma_{i,j}^d\cup \Gamma_{i,j}^u$ denote the union of
the following sets of ``paths'' from $(1,j)$ to $(i,1)$:

\begin{eqnarray}
&\Gamma_{i,j}^d=\{[i_1,\dots,i_r;j_1,\dots,j_r]\mid r\in {\mathbb
  N},\\\nonumber & 1= i_1<i_2<\dots<i_r= i, 1\leq j_r<\dots <j_1=j\},\\\nonumber 
&\Gamma_{i,j}^u=\{[i_1,\dots,i_r;j_1,\dots,j_r]\mid r\in {\mathbb
  N},\\\nonumber & 1< i_1<i_2<\dots<i_r= i, 1\leq j_r<\dots <j_1= j\}. 
\end{eqnarray}

For $j=1, i\geq 1$ we interpret the above as
$\Gamma_{i,1}^d=\emptyset$ and $\Gamma_{i,1}^u=\{[i;1]\}$ whereas for
$i=1, j>1$ it is $\Gamma_{1,j}^u=\emptyset$ and
$\Gamma_{1,j}^u=\{[1;j]\}$.

\smallskip

For each $z_{i,j}$ we define an operator ${\s K}_{i,j}$ in analogy
with (\ref{43}) and an operator $K_{i,j}$ in analogy with (\ref{44}),
and finally we set ${\s O}_{i,j}={\s K}_{i,j}K_{i,j}^2$.

For $g=[1_1,\dots,i_r;j_1,\dots, j_r]\in \Gamma^d_{i,j}$, set
\begin{eqnarray}\nonumber S_g&=&\{(s,t)\mid \exists x=1,\dots,   r-1: s=i_x\text{ and }
  j_{x+1}<t<j_x\},\\\nonumber  T_g&=&\{(s,t)\mid \exists x=1,\dots,  r-1:
  t=j_{x+1}\text{ and }  i_{x}<t<i_{x+1}\}, \\ \text{and}\\\nonumber
  D^d_{i,j}(g)&=&\left(\prod_{y=1}^{r-1}{\s O}_{i_y,j_{y+1}}\prod_{(s,t)\in
    S_g}K_{s,t}\prod_{(s,t)\in T_g}K_{s,t}^{-1}\prod_{x<j_r}
  K_{i,x}^{-1}\right)\prod_{x=1}^r\frac{\partial}{\partial
  z_{i_x,j_x}}. 
\end{eqnarray}  

Likewise, for $g=[1_1,\dots,i_r;j_1,\dots, j_r]\in \Gamma^u_{i,j}$,
set, for convenience, $i_0=0$ and
\begin{eqnarray}\nonumber U_g&=&\{(s,t)\mid \exists x=1,\dots,   r-1: s=i_{x+1}\text{ and }
  j_{x+1}<t<j_x\},\\\nonumber  V_g&=&\{(s,t)\mid \exists x=1,\dots,  r:
  t=j_{x}\text{ and }  i_{x-1}<t<i_{x}\}, \\ \text{and}\\\nonumber
  D^u_{i,j}(g)&=&\left( \prod_{y=1}^{r-1}{\s O}_{i_y,j_{y+1}}\prod_{(s,t)\in U_g}K_{s,t}\prod_{(s,t)\in V_g}K_{s,t}^{-1}\prod_{x<j_r}K_{i,x}^{-1}\right)\prod_{x=1}^r\frac{\partial}{\partial
  z_{i_x,j_x}}. 
\end{eqnarray}  

Then

\begin{equation}
\frac{\partial}{\partial X_{i,j}}=\sum_{g\in \Gamma_{i,j}^d}D^d_{i,j}(g)+\sum_{g\in \Gamma_{i,j}^u}D^u_{i,j}(g).
\end{equation}

Observe that the lowest order differential operator occurring as a
summand in $\frac{\partial}{\partial X_{i,j}}$ is
$\left(\prod_{y<j_r}K_{1,y}^{-1}\prod_{x<i}K_{x,j}^{-1}\right)\frac{\partial}{\partial
  z_{i,j}}$.

\medskip

It is not clear if an analogue of the operator $G$ exists for
higher order algebras.

\medskip

\bibliographystyle{alpha} \bibliography{qt23-nyeste}

\end{document}